\documentclass[12pt,a4paper]{amsart}
\usepackage{amsthm}
\usepackage{amsmath}
\usepackage{amsfonts}
\usepackage{amssymb}
\usepackage[left=2.5cm, right=2.5cm, top=3.0cm, bottom=3.0cm]{geometry}
\usepackage{color}
\usepackage{bbm}
\usepackage[arrow, matrix, curve]{xy}
\usepackage{graphicx}
\usepackage[colorlinks=true,linkcolor=green]{hyperref}
\makeatletter
\renewcommand*{\eqref}[1]{
	\hyperref[{#1}]{\textup{\tagform@{\ref*{#1}}}}}
\makeatother

\newcommand{\dx}{\dot{x}}
\newcommand{\ddx}{\ddot{x}}

\newcommand{\beq}{\begin{equation}}
	\newcommand{\beqn}{\begin{equation*}}

		\newcommand{\dgamma}{\dot{\gamma}}

		\newcommand{\eeq}{\end{equation}}
	\newcommand{\eeqn}{\end{equation*}}
	\newcommand{\grad}{\textrm{grad}}
 \newcommand{\kin}{T}

\newcommand{\R}{\mathbb{R}}

\newcommand{\Stau}{S^1_{\tau}}
\newcommand{\tgamma}{\tilde{\gamma}}

\newcommand{\tx}{\widetilde{x}}

\theoremstyle{plain}
\newtheorem{theorem}{Theorem}
\newtheorem{proposition}{Proposition}
\newtheorem{corollary}{Corollary}
\newtheorem{lemma}{Lemma}

\theoremstyle{definition}

\theoremstyle{remark}
\newtheorem{remark}{Remark}
\newtheorem{example}{Example}
\newtheorem{number-env}[theorem]{}
\title[Periodic orbits]{Periodic 
	orbits of 
	reversible Lagrangian systems without self-intersections
	and Ma\~n\'e genericity}
\author{Hans-Bert Rademacher}
\address{Mathematisches Institut, 
	Universit{\"a}t Leipzig, D--04081 Leipzig, Germany}
\email{rademacher@math.uni-leipzig.de}
\urladdr{\url{www.math.uni-leipzig.de/~rademacher}}
\date{2026-02-19
}
\begin{document}
	\begin{abstract}
 Bernard~\cite{Be}  showed that a
		Ma\~n\'e generic convex Hamiltonian has
		only non-degenerate periodic orbits on a given
		energy level.
		We show that one can use this result to prove that
		for a generic potential the prime periodic orbits
        of fixed energy of a Lagrangian system of classical
        type on a compact manifold of dimension $n\ge 3$  do not
        have self-intersections and do not intersect each other.       
	\end{abstract}
	\maketitle
	\baselineskip 18pt
	\section{Introduction}

 We consider a \emph{reversible Lagrangian system of classical type}
 resp. a \emph{natural mechanical system} defined by
 a \emph{reversible Finsler metric} $F:TM \longrightarrow \R$ on the tangent bundle
 $TM$ of a differentiable manifold $M$ of dimension
 $n$ with foot point map $\tau:TM \longrightarrow M,$
 and a \emph{potential energy} $U: M\longrightarrow \R,$
 i.e. a smooth function. 
 Then $\kin: TM \longrightarrow \R, \kin(v)=F^2(v)/2$ 
 is the \emph{kinetic energy.}
 The manifold represents the 
 configuration space of a dynamical system, the
 \emph{Lagrangian} $L: TM \longrightarrow \R$ is defined by
 \begin{equation*}
     v \in TM \longmapsto L(v)=
     \kin(v)-U(\tau(v))=\frac{1}{2} F^2(v)-U(\tau(v))\in \R\,.
 \end{equation*} 
 The orbits  $t \mapsto \dgamma(t)\in TM$
 of the Lagrangian system are solutions of 
 the Lagrange equations~\eqref{eq:lagrange}.
 The \emph{total energy} is given by
 \begin{equation*}
  \tilde{H}: TM \longrightarrow \R;
 \tilde{H}(v)=\kin(v)+U(\tau(v))	\,.
\end{equation*}
 For an orbit $\dgamma=\dgamma(t)$ the total
 energy is constant, i.e. there is $E \in \R$
 with
 \begin{equation*}
     E=
     \tilde{H}(\dgamma(t))=\kin(\dgamma(t))+U(\gamma(t))=\frac{1}{2} F^2(\dgamma(t))+U(\gamma(t))
 \end{equation*}
 for all $t.$  
 
If $E$ is a regular value
 of the potential $U$ then the projection $t \mapsto \gamma(t)=\tau(\dgamma(t))\in M$ of 
 the orbit $t \mapsto \dgamma (t)$ of 
 the Lagrangian system is (up to parametrization) a
 geodesic curve of the 
 Finsler metric
 \begin{equation}
 	\label{eq:Jacobi-Finsler}
 \overline{F}^2(v)=
 2\left\{E-U(\tau(v))\right\}\, F^2(v)
 \end{equation}
 defined on the open subset 
 $M^{<E}=U^{-1}(-\infty,E)$ of $M$ 
 called \emph{potential well} or \emph{Hill's region}
 with smooth
 boundary $M^E=\{x\in M\,;\, U(x)=E\}=\partial M^{<E}.$
 Let $M^{\le E}:=\{x\in M\,;\, U(x)\le E\}=
 M^{<E}\cup M^E.$ 
 The metric $\overline{F}$ is called 
 \emph{Jacobi(-Maupertuis) Finsler metric,}
 this result is well known in particular if $F$ is a Riemannian metric,
 cf. for example \cite[Sec. 4.1]{AKN}. 
 For completness
 a detailed computation valid also for a Finsler metric
 is presented in the following 
 Section~\ref{sec:finsler}
 and Proposition~\ref{pro:jacobi-finsler}.  
 Note that the metric
 $\overline{F}$ is conformally equivalent to 
 $F$ and degenerates
 on $M^E$ if $E< \max U.$ 
 Using the \emph{Legendre transformation}
 \begin{equation}
 	\label{eq:Legendre}
    \mathcal{L}: TM \longrightarrow T^*M, v \mapsto \mathcal{L}(v);
    \mathcal{L}(v)(w)=d F^2_v(w)
    =\left.\frac{d}{dt}\right|_{t=0}F^2(v+tw)
 \end{equation}
between the tangent bundle $TM$ and the cotangent bundle $T^*M$
 there is an associated \emph{Hamiltonian system of
 classical type} given by the 
 \emph{Hamiltonian} $H:T^*M \longrightarrow \R$
 defined by $H(\omega)=\kin(\mathcal{L}^{-1}(\omega))+U(x),
 \omega \in T_x^*M,$
 i.e. $\tilde{H}(v)=H(\mathcal{L}(v)).$
 Under the Legendre transformation the
 orbits of the Hamiltonian system correspond to the orbits
 of the Lagrangian system.
 With a periodic orbit 
 $t\in \Stau=\R/(\mathbb{Z}\tau) \mapsto \dgamma(t)\in TM$
 of minimal period $\tau>0$
 for $m \in \mathbb{Z}$
 the iterates $t \in \Stau \mapsto \dgamma (mt)$
 are periodic orbits, too. We call a periodic orbit
 of minimal period $\tau$
 \emph{prime,} i.e. the mapping
 $t \in \R/(\mathbb{Z}\tau)\mapsto \dgamma (t)\in TM$ is injective.
 In particular  
 it is not the iterate for some $m>1$
 of another periodic orbit.
 We call a periodic orbit $t \in \Stau \mapsto \dgamma (t)
 \in TM$ \emph{simple} or
 \emph{without self-intersection} if
 for all $s,t \in \Stau, s\not=t$ with $\gamma(s)=\gamma(t):
 \dgamma(t)=-\dgamma(s).$  
 Hence a simple periodic orbit is prime.
 For a periodic orbit
 $\dgamma: \Stau \longrightarrow TM$ of period $T>0$  
 of the reversible Lagrangian system
 of classical type of total enery $E,$ which is a regular
 value of the potential $U,$ there are
 two cases, cf. \cite[Prop.4.5]{AKN}:
 \begin{itemize}
 	\item[(a)] 
 	If $\dgamma (t)\not=0$ for all $t\in \Stau,$ then
 	the projection
 	$t \mapsto \gamma(t)=\tau(\dgamma(t))$ of a
 	 periodic orbit $t \mapsto \dgamma(t)$ 
 	 of the Lagrange system corresponds to
 	(up to parametrization) a periodic geodesic
 	(also called \emph{closed geodesics})
 	of the Jacobi-Finsler metric $\overline{F}$
 	defined by 
 	Equation~\eqref{eq:Jacobi-Finsler} on $M.$
 	These periodic orbits are called
 	\emph{rotations} in \cite[Sec.4.2.1]{AKN}.
 	For a rotation 
 	$t \in \Stau \mapsto \dgamma(t)\in TM$ of 
 	minimal period $T$
 	we have the following property:
 	$\gamma(s)=\gamma(t)$ for $s\not=t$
 	implies that $\dgamma(s),\dgamma(t)$ are linearly independent.
 	Hence for a simple rotation 
 	$t\in \Stau\mapsto \dgamma (t)\in TM$ the 
 	map $t\in \Stau \mapsto \gamma (t)\in M$ is injective.
 	\item[(b)] 
 	 	If there is $t_0 \in \Stau$ with 
 	$\dgamma(t_0)=0,$ i.e. $U(\gamma(t_0))=E,$ 
 	then $\gamma$ is a periodic orbit with \emph{rest point}
 	$\gamma(t_0).$ 
 	If follows that $E< \max U,$
 	the orbit $\gamma$ lies in the potential well
 	$M^{\le E}.$
 	For simplicity we assume 
 	$t_0=0,$ i.e. $\dgamma(0)=0.$ 	
    Since the Lagrangian resp.
 	the Finsler metric $F$ is reversible it follows that
 	$\gamma(\tau/2+t)=\gamma(\tau/2-t)$ for all $t.$
 	Such a periodic orbit is called
 	\emph{brake orbit.} 
 	Hence a brake orbit oscillates between
 	the rest points $\gamma(0)$ and $\gamma(\tau/2)$ with
 	$E=U(\gamma(0))=U(\gamma(\tau/2))$ 
 	and $\dgamma(0)=0, \dgamma(\tau/2)=0.$
 	A brake orbit $\gamma$ is prime, if the period $\tau$ is minimal,
 	i.e. there are no further rest points except
 	$\gamma(0), \gamma(\tau/2).$ 
 	
 	A brake orbit $t \in \Stau\longmapsto \dgamma(t)\in TM$
 	of period $\tau$ is
 	\emph{simple,} or \emph{without self-intersections,} if 
 	$\gamma(s)=\gamma(t), s\not=t, s,t \in \Stau$ implies 
  $\dgamma(s)=-\dgamma(t),$ resp. $s+t=\tau.$
 	Hence in this case the
 	restriction $\gamma: [0,\tau/2]\longrightarrow M$
 	is injective, and $\gamma: [0,\tau]\longrightarrow M$
  satisfies $\gamma(\tau/2-t)=\gamma(\tau/2+t)$
  for all $t,$
  i.e. $\gamma: (0,\tau) \longrightarrow M$
  is a $(2:1)$-map, where $\gamma(s)=\gamma(t), t\not=s$
  implies $\dgamma(s)=-\dgamma(t).$
  Brake orbits are also called \emph{librations,}
  cf. for example~\cite{BK}, or \emph{round trip orbits,}
  cf.~\cite[Prop.1]{Be}.
   \end{itemize}
 After preliminary results by Oliveira~\cite{Ol} and
 Rifford and Ruggiero~\cite[Theorem 1.2]{RR} Bernard 
 proved in~\cite[Theorem 1]{Be}
 that for a Ma\~n\'e generic convex
 Hamiltonian periodic orbits on a given regular
 energy level are non-degenerate.
 Using this result we obtain the following main result of
 this paper:
 \begin{theorem}
 	\label{thm:main}
 	Let $F$ be a reversible Finsler metric on the compact smooth
 	manifold $M$ 
 	of dimension $n=\dim M\ge 3$
 	defining the kinetic energy 
 	$\kin=F^2/2$ and let $E \in \R.$
 	Then for a generic function $U \in C^{\infty}(M)$
 	the number
 	$E$ is a regular value of the function $U$ as well as of
 	the total energy $\tilde{H}=\kin+U,$
 	and the following holds:
 	
 	All periodic orbits of the reversible Lagrangian system 
 	of classical type defined
 	by the Lagrangian $L=\kin-U$ of total energy $E$ 
 	are non-degenerate and all prime
 	periodic orbits of total energy $E$ do not have
 	self-intersections and they do not intersect each other. 	
 \end{theorem}
Hence for a Ma\~n\'e generic reversible Lagrange system
of classical type defined on the tangent bundle of a compact
manifold of dimension $\ge 3,$
i.e. for generic potential,
the prime periodic orbits of 
fixed energy are simple and they do not intersect each
other.
Using the Legendre transformation one can also state
a Hamiltonian version of this result.

For a reversible Lagrangian system of classical type 
on a compact manifold $M$
with regular value $E$
there exists a periodic motion,
this was shown by Gluck and Ziller~\cite[Thm.A]{GZ},
cf. also~\cite[Thm.1]{Kn}, \cite{Ko} and \cite[Sec.4.2]{AKN}.
If $M^E$ is non-empty then a brake orbit exists.
If $M^E$ has 
$k$ connected components it was shown by Bolotin and Kozlov~\cite{BK}
that there are $k$ geometrically distinct 
brake orbits, which do not have self-intersections.

As a consequence of
Theorem~\ref{thm:main} we obtain:
\begin{corollary}
	\label{cor:periodic}
	Let $M$ be a compact differentiable manifold of dimension
	$n\ge 3$ with a reversible Finsler metric $F$ defining the
	kinetic energy $\kin=F^2/2$ and let $E \in \R.$
 	Then there is an open and dense subset $\mathcal{U}$ of
 the set $C^{\infty}(M)$ of smooth functions such that the following
 holds:

The number $E$ is a regular value of $U \in \mathcal{U},$ and
a regular value of the total energy $\tilde{H}=\kin+U,$
and there exists a non-degenerate periodic orbit 
 of energy $E$ of
 the Lagrangian system defined by  
 the Lagrangian $L=\kin-U,$
 which is simple. 
 
 Furthermore, if $M^E$ is non-empty
 then there exists a brake orbit of energy $E,$ which is
 non-degenerate and simple.
 \end{corollary}
\begin{remark} 
\label{rem:seifert}
Seifert proves in~\cite{Se} the existence of a brake orbit
if the potential well is diffeomorphic to an $n$-disc,
compare also~\cite{GZ} and
\cite{We}. Seifert also states in a footnote~\cite[p.15]{Se}
that it might be possible to prove the existence of
$n$ geometrically distinct brake orbits. 
Under additional assumptions this result has been shown,
cf. for example~\cite[Theorem D]{GZ}
or \cite{Ko}.
The result without additional assumption is announced by
Giambo, Giannoni and Piccione~\cite{GGP2022}. 
We can use Theorem~\ref{thm:main}
and obtain that for an open and dense set of potentials
in these results
the $n$ geometrically distinct prime brake orbits are
in addition non-degenerate and simple and do not
intersect each other.
\end{remark}
The starting point for this investigation was a question 
by Alin Albu-Sch\"affer about self-intersections of periodic
orbits. 
Remark~\ref{rem:seifert} gives a partial answer to this question.
I am grateful to Alin Albu-Sch\"affer for 
several interesting discussions.
For a connection between brake orbits with
robotics planning see for example~\cite{AS}.
\begin{remark}
(a) The author has shown in~\cite[Theorem 1]{R24} that on 
	a compact manifold $M$ of dimension $n\ge 3$ for a generic
	Riemannian metric or for a generic
	reversible Finsler metric the prime
	closed geodesics are simple and do not intersect each other.
	Analogous results for generic Riemannian metrics and
	geodesic loops resp. orthogonal geodesic chords in 
	manifolds with boundary are shown in~\cite{R24a}.
	In this paper we perturb the potential which corresponds
	to a conformal change of the corresponding Jacobi-Finsler metric.
	This approach goes back to Ma\~n\'e, cf.~\cite{Ma}.

 \medskip
	 
	(b) Bernard's result~\cite[Theorem 1]{Be} also shows that a bumpy 
	metrics theorem holds in the restricted class of 
	conformally equivalent metrics:
	
	Let $M$ be a compact differentiable manifold
		of dimension $n\ge 2$ with Riemannian metric $g,$ resp.
		reversible Finsler metric $F.$ Then in an arbitrary neighborhood
		of the zero function in the space $C^{\infty}(M)$ there 
		is a smooth function $\phi$ such that the conformally
		equivalent metric $\overline{g}=\exp(\phi) g$ is bumpy,
		i.e. all closed geodesics are non-degenerate.		

  \medskip	

	(c) The reversibility of the Finsler metric is
	essential for our proof since for non-reversible Finsler metrics the statements
	about double points of periodic orbits obtained
	in Lemma~\ref{lem:intersection} do not hold. 

 \medskip

 (d)
 Recently Dahinden and de Pooter show in~\cite[Thm.1, Cor.3]{DdP}
 that also for a generic \emph{non-reversible} Finsler metric on a compact manifold
 the closed geodesics do not have self-intersections.
	
\end{remark}
\section{Finsler metrics, Lagrangian and Hamiltonian systems of classical type}
\label{sec:finsler}
On a smooth manifold $M$ we choose
a local coordinate
system $x=(x^1,\ldots,x^n)$ which then
induces an associated local coordinate
system
$(x,\dx)=(x^1,\ldots,x^n,\dx^1,\ldots,\dx^n)$ on the tangent
bundle $TM.$ Hence the coordinate fields
$\partial_i=\partial/\partial x^i,i=1,\ldots,n$
define a basis of the tangent space $T_xM.$
Then a \emph{Finsler metric}
$F:TM \longrightarrow \R$ is a continuous function,
which is smooth outside the zero section,
positively homogeneous, i.e. $F(x,\lambda\dx)=
\lambda F(x,\dx)$
for all tangent vectors $(x,\dx)$
and $\lambda>0,$
and satisfies the \emph{Legendre condition:}
For all $(x,\dx)$ with $\dx\not=0$ the matrix defined by
\begin{equation}
\label{eq:gij}		
g_{ij}(x,\dx)=
\frac{1}{2} 
\frac{\partial^2 F^2(x,\dx)}{\partial \dx^i \partial \dx^j}	,
1\le i,j \le n
\end{equation}
is positive definite. Since $F$ is positively homogeneous
the function $(x,\dx)\mapsto g_{ij}(x,\dx)$
is positively homogenous of degree $0$
with respect to $\dx.$
The Finsler metric is called \emph{reversible} if
$F(x,-\dx)=F(x,\dx).$
The \emph{Cartan tensor} $C_{ijk}=C_{ijk}(x,\dx)$
defined by
\begin{equation*}
C_{ijk}(x,\dx)=
\frac{1}{2}\frac{\partial g_{ij}(x,\dx)}{\partial \dx^k}	=
\frac{1}{4} 
\frac{\partial^3 F^2(x,\dx)}{\partial\dx^i\partial\dx^j\partial\dx^k}
\end{equation*} 
is symmetric in $i,j,k$ and homogenous of degree $-1$ in $\dx,$ hence
\begin{equation}
	\label{eq:cartan-symmetry}
C_{ijk}(x,\dx)\dx^i=C_{ijk}(x,\dx)\dx^j=C_{ijk}(x,\dx)\dx^k=0\,,
\end{equation}
cf. \cite[p.84]{Sh}.
Here we use the Einstein summation convention.
The Cartan tensor vanishes identically if and only if
the Finsler metric is a Riemannian metric.

The \emph{Christoffel symbols 
$\gamma_{ijl}$ of the first kind}
are defined by
\begin{equation*}
\gamma_{ijl}(x,\dx)=\frac{1}{2}
\left\{
\frac{\partial g_{li}(x,\dx)}{\partial x^j}+
\frac{\partial g_{jl}(x,\dx)}{\partial x^i}-
\frac{\partial g_{ij}(x.\dx)}{\partial x^l}
\right\}
\end{equation*}
Since 
$\dx^j \dx^k$ resp.
$\dx^i\dx^j\dx^k$ is an expression
symmetric in 
$j,k$ resp.
$i,j,k$ we obtain the following identities,
which we will use later:
\begin{equation}
	\label{eq:christoffel}
	\gamma_{ijk}(x,\dx)\dx^j\dx^k=
	\frac{1}{2}
	\left\{
	\frac{\partial g_{ki}(x,\dx)}{\partial x^j}
	+\frac{\partial g_{jk}(x,\dx)}{\partial x^i}
	-\frac{\partial g_{ij}(x,\dx)}{\partial x^k}	
	\right\} \dx^j \dx^k
	=\frac{\partial g_{jk}(x,\dx)}{\partial x^i}
		\dx^j \dx^k\,,		
\end{equation}
and
\begin{equation}
\gamma_{ijk}(x,\dx)\dx^i\dx^j \dx^k
=\frac{1}{2}\frac{\partial g_{ij}(x,\dx)}{\partial x^k}
\dx^i \dx^j \dx^k\,.
\end{equation}
The \emph{Christoffel symbols 
$\Gamma_{ij}^k$ of the second kind} are defined
by
\begin{equation}
\Gamma_{ij}^k(x,\dx)=
g^{kl}(x,\dx)\gamma_{ijl}(x,\dx)=
\frac{1}{2}g^{kl}(x,\dx)\left\{
\frac{\partial g_{li}(x,\dx)}{\partial x^j}+
\frac{\partial g_{jl}(x,\dx)}{\partial x^i}-
\frac{\partial g_{ij}(x.\dx)}{\partial x^l}
\right\}\,.
\end{equation}
The \emph{geodesic coefficients} 
$G^k(x,\dx),k=1,\ldots,n\,$ are
defined by
\begin{equation}
\nonumber
G^k(x,\dx)=\frac{1}{4}g^{kl}(x,\dx)\left\{2 
\frac{\partial g_{li}(x,\dx)}{\partial x^j}-
\frac{\partial g_{ij}(x,\dx)}{\partial x^l}
\right\}
\dx^i \dx^j
=\frac{1}{2}\Gamma_{ij}^k(x,\dx)\dx^i\dx ^j\,.
\end{equation}
Then the smooth curve $t\mapsto \gamma(t)=(x^1(t),\ldots, x^n(t))\in M$ 
is a \emph{geodesic (curve)} of the Finsler metric $F$
if and only if 
$$
\ddx^k=\ddx^k(t)
=\frac{d^2}{dt^2}x^k(t)=\frac{d}{dt}\dx^k (t)
$$
satisfies 
\begin{equation}
\label{eq:geodesic}
\ddx^k(t)+2 G^k(x(t),\dx(t))=
\ddx^k(t)+\Gamma^k_{ij}(x(t),\dx(t))\dx^i(t)\dx^j(t)=0
\end{equation}
resp.
\begin{equation*}
g_{ik}(x(t),\dx(t))\ddx^k(t)
+\gamma_{ijl}(x(t),\dx(t))\dx^j(t)\dx^l(t)=0\,.    
\end{equation*}
It follows that a geodesic $x=x(t)$ is parametrized
proportional to arc length, as the following
computation shows:
\begin{eqnarray*}    
\frac{d}{dt}F^2(x(t),\dx(t))=
\frac{1}{2}\frac{d}{dt}\left\{
g_{ij}(x(t),\dx(t))\dx^i(t)\dx^j(t)
\right\}
\\
=\frac{1}{2}
\frac{\partial g_{ij}(x,\dx)}{\partial x^k}
\dx^k \dx^i \dx^j
+
\frac{1}{2}
\frac{\partial g_{ij}(x,\dx)}{\partial \dx^k}
\ddx^k \dx^i \dx^j
+g_{ij}(x,\dx)\ddx^i \dx^j
\\
= \gamma_{ijk}\dx^i\dx^j\dx^k+g_{ij}(x,\dx)\left(
-\Gamma_{il}^k(x,\dx)\dx^i\dx^j 
\right)
\dx^l=0\,.
\end{eqnarray*}
Hence $F^2(x(t),\dx(t))=F^2(x(0),\dx(0))$ for all $t.$
A \emph{classical reversible mechanical system} on a differentiable
manifold $M$ resp.
a \emph{reversible Lagrangian system of classical type} is given by
a reversible Finsler metric $F$ on $M$
defining the 
\emph{kinetic energy} $
\kin(x,\dx)=F^2(x,\dx)/2$
and a smooth function,
called \emph{potential}
$U:M \longrightarrow \R.$
The \emph{Lagrangian} $L:TM \longrightarrow \R$
is then defined by
\begin{equation}
	L(x,\dx)=\kin(x,\dx)-U(x)=\frac{1}{2} F^2(x,\dx)-U(x)=
	\frac{1}{2}g_{ij}(x,\dx)\dx^i\dx^j-U(x)\,.	
\end{equation}
where $g_{ij}(x,\dx)$ are defined by Equation~\eqref{eq:gij}.
\emph{Solutions} $x=x(t)$
resp. \emph{orbits} 
$(x,\dx)=(x(t),\dx(t))=(x(t),\dx(t))$ of the Lagrangian
system are solutions of the 
\emph{Lagrange equations:}
\begin{equation}
\label{eq:lagrange}
\frac{d}{dt}\frac{\partial L(x,\dx)}{\partial \dx^i}=
\frac{\partial L(x,\dx)}{\partial x^i},
i=1,\ldots,n\,.
\end{equation}
For the Lagrangian $L(x,\dx)=\kin(x,\dx)-U(x)$ of a classical system
we obtain from Equation~\eqref{eq:lagrange}:
\begin{eqnarray*}
\frac{d}{dt}\frac{\partial L(x,\dx)}{\partial \dx^i}=
\frac{1}{2}\frac{d}{dt}\frac{\partial F^2(x,\dx)}{\partial \dx^i}
=\frac{d}{dt} \left(g_{ik}(x,\dx)\dx^k\right)\\
=\frac{\partial g_{ik}(x,\dx)}{\partial x^l}\dx^l \dx^k
+\frac{\partial g_{ik}(x,\dx)}{\partial \dx^l}\ddx^l\dx^k	
+g_{ik}(x,\dx)\ddx^k
=\frac{\partial g_{ik}(x,\dx)}{\partial x^l}\dx^l \dx^k
+g_{ik}(x,\dx)\ddx^k\\=
\frac{\partial L(x,\dx)}{\partial x^i}
=\frac{1}{2}\frac{\partial F^2(x,\dx)}{\partial x^i}-
\frac{\partial U(x)}{\partial x^i}
=\frac{1}{2}
\frac{\partial g_{kl}(x,\dx)}{\partial x^i}\dx^k\dx^l
-\frac{\partial U(x)}{\partial x^i}\,.
\end{eqnarray*}
Hence we obtain the following differential equation
of second order for an orbit $x=x(t)$ of the Lagrangian
system, generalizing the geodesic equation~\eqref{eq:geodesic}:
\begin{eqnarray}
\nonumber
\ddx^i+\frac{1}{2}g^{ik}(x,\dx)
\left\{
2\frac{\partial g_{jk}(x,\dx)}{\partial x^l}-
\frac{\partial g_{jl}(x,\dx)}{\partial x^k}
\right\}\dx^j\dx^l
+g^{ik}(x,\dx)\frac{\partial U(x)}{\partial x^k}
\\
\label{eq:orbit-lagrange}
=\ddx^i +2 G^i(x,\dx)+
g^{ik}(x,\dx)\frac{\partial U(x)}{\partial x^k}\\
\label{eq:orbit-lagrange1}
=\ddx^i +\Gamma^i_{jk}(x,\dx)(x,\dx)+
g^{ik}(x,\dx)\frac{\partial U(x)}{\partial x^k}=0\,,
\end{eqnarray}
resp.
\begin{equation}
	\label{eq:orbit-lagrange3}
	g_{ik}(x,\dx)\ddx^i +\gamma_{ijk}(x,\dx)(x,\dx)+
	\frac{\partial U(x)}{\partial x^i}=0\,.
\end{equation}
Along a solution $x=x(t)$ of a Lagrangian system the total
energy $T(x,\dx)+U(x)$
is constant, which follows from
the following computation:
\begin{eqnarray*}
\frac{d}{dt}
\left\{
\frac{1}{2} F^2(x(t),\dx(t))+U(x(t))
\right\}
=
\frac{d}{dt} 
\left\{\frac{1}{2} g_{ij}(x,\dx) \dx^i\dx^j +U(x(t))\right\}\\
=
\frac{1}{2} \frac{\partial g_{ij}(x,\dx)}{\partial x^k}
\dx^i\dx^j\dx^k
+\frac{\partial g_{ij}(x,\dx)}{\partial \dx^k}
\ddx^k \dx^i\dx^j
+ g_{ij}(x,\dx)\ddx^i\dx^j+
\frac{\partial U}{\partial x^k}\dx^k\\
=\frac{1}{2} \frac{\partial g_{ij}(x,\dx)}{\partial x^k}
\dx^i\dx^j\dx^k
+ g_{ij}(x,\dx)\ddx^i\dx^j+
\frac{\partial U}{\partial x^k}\dx^k=0\\
\end{eqnarray*}
Here in the last line we use
Equation~\eqref{eq:orbit-lagrange}
and Equation~\eqref{eq:cartan-symmetry}.

Hence for a solution
$x=x(t)$ 
of the Lagrangian system of total energy
$E \in \R$ we obtain
\begin{equation}
	\label{eq:energy-orbit}
	\kin(x(t),\dx(t))=\frac{1}{2}F^2(x(t),\dx(t))=E-U(x(t))\,.
\end{equation}
The solutions of
a Lagrangian system of classical type can be seen as geodesics
of a Finsler metric $\overline{F}$
in the conformal class of the Finsler
metric $F$ of the Lagrangian system (up to parametrization),
cf. Proposition~\ref{pro:jacobi-finsler}.
This metric is also
called \emph{Jacobi-Finsler metric,}
it is given by Equation~\eqref{eq:Jacobi-Finsler}.
The \emph{Legendre transformation} 
$$\mathcal{L}:(x,\dx)=
(x^1,\ldots,x^n,\dx^1,\ldots,\dx^n) \in TM \longmapsto (x,y)=(x^1,\ldots,x^n,y_1,\ldots,y_n)\in T^*M$$
with respect to the coordinate fields $\partial_1,\ldots,\partial_n$
of $TM$ resp. its dual basis
$dx_1,\ldots,dx_n$ of $T^*M$ is given by:
\begin{equation*}
y_i=\frac{\partial{L(x,\dx)}}{\partial \dx^i}
=
g_{ij}(x,\dx)\dx^j
+\frac{\partial g_{ij}(x,\dx)}{\partial \dx^l}\dx^l
= g_{ij}(x,\dx)\dx^j,
\end{equation*}
here we use again Equation~\eqref{eq:cartan-symmetry}.
Then we obtain for the total energy 
\begin{eqnarray*}
H(x,y)=\frac{1}{2}F^2(x,\dx) +U(x)=
\frac{1}{2}g_{ij}(x,\dx)\dx^i\dx^j+U(x)\\
=y_i x^i+U(x)=\frac{1}{2}g^{ij}(x,y) y_i y_j +U(x)\,.
\end{eqnarray*}
\begin{remark}
On a Riemannian manifold with Levi-Civita connection $\nabla$
and covariant derivative $\frac{\nabla}{dt}$ along a curve $t \longmapsto \gamma(t) \in M$
Equation~\eqref{eq:orbit-lagrange}
has the following invariant form:
$$
\frac{\nabla}{dt}\dgamma(t)+\grad U(\gamma(t))=0\,.
$$
Here $\grad U$ is the \emph{gradient} of $U.$
\end{remark}
A \emph{Hamiltonian system of classical type}
is given by a Finsler metric $F^*$ on the cotangent bundle
$T^*M$ of a differentiable manifold $M$ defining the kinetic
energy and a smooth function $U,$
called the potential. Hamiltonian systems of classical
type correspond via the Legendre transformation to
Lagrangian systems of classical type.

Then we obtain for a solution of the Lagrangian system resp.
the Hamiltonian system
using Equation~\eqref{eq:cartan-symmetry},
Equation~\eqref{eq:orbit-lagrange3},
Equation~\eqref{eq:christoffel}:
\begin{eqnarray*}
\dot{y}_i=\frac{d}{dt}\left(
g_{ij}(x,\dx)\dx^j\right)=
\frac{\partial g_{ij}(x,\dx)}{\partial x^k}
\dx^j\dx^k+ 
\frac{\partial g_{ij}(x,\dx)}{\partial \dx^k}\dx^j\ddx^k
+g_{ij}(x,\dx)\ddx^j+g_{ij}(x,\dx)\ddx^j=\\
=\frac{\partial g_{ij}(x,\dx)}{\partial x^k}
\dx^j\dx^k
-\gamma_{ijk}(x,\dx)\dx^j\dx^k
+\frac{\partial U}{\partial x^i}=
\frac{1}{2}\frac{\partial g_{jk}(x,\dx)}{\partial x^i}
\dx^j\dx^k
+\frac{\partial U}{\partial x^i}
=\frac{\partial H(x,\dx)}{\partial x^i}
=\frac{\partial H(x,y)}{\partial x^i}
\,; 
\end{eqnarray*}
\begin{eqnarray*}
\frac{\partial H(x,y)}{\partial y_i}	
=\frac{\partial }{\partial y_i}
\left\{\frac{1}{2}g^{kj}(x,y)y_ky_j
+U(x)
\right\}\\
=
\frac{1}{2}\frac{\partial g^{kj}(x,y)}{\partial y_i}y_ky_j+
g^{ij}(x,y)y_j
=g^{ij}(x,y)y_j=\dx^i,
\end{eqnarray*}
here we use that $\partial{g^{kj}(x,y)}/{\partial y_i}$ 
is homogeneous of
order $-1$ with respect to $y$ which implies
\begin{equation*}
	\frac{\partial g^{kj}(x,y)}{\partial y_i}y_ky_l=0\,,	
\end{equation*}
cf. Equation~\eqref{eq:cartan-symmetry}.
Hence we obtain \emph{Hamilton's equations:}
\begin{equation*}
\dx^i=\frac{\partial H(x,y)}{\partial y_i}	\,;\,	
\dot{y}_i=-\frac{\partial H(x,y)}{\partial x^i}	\,;\,
i=1,\ldots,n\,.
\end{equation*}
Let $\phi: M\longrightarrow \R$ be a smooth function.
The geodesic coefficients $\overline{G}^i(x,\dx)$ of the
conformally equivalent Finsler metric
$\overline{F}(x,\dx)=\exp(\phi(x)) F(x,\dx)$
are given by:
\begin{eqnarray}
\overline{G}^i(x,\dx)=
\frac{1}{4}
\exp(-\phi(x))g^{ik}(x,\dx)
\left\{2\frac{\partial(\exp(\phi(x))g_{jk}(x,\dx))}{\partial x^l}
\right.
\nonumber
\\
- \left.\frac{\partial(\exp(\phi(x))g_{jl}(x,\dx))}{\partial x^k}
\right\}\dx^j \dx^l \nonumber\\
=
G^i(x,y)+\frac{1}{4}
g^{ik}(x,\dx)
\left\{2\frac{\partial \phi}{\partial x^l}g_{jk}(x,\dx)-
\frac{\partial \phi}{\partial x^k} g_{jl}(x,\dx)\right\}	
\dx^j \dx^l
\nonumber
\\
=
G^i(x,\dx)+
\frac{1}{4}
\left\{2
\frac{\partial \phi(x)}{\partial x^l}\dx^l\dx^i
-g^{ik}(x,\dx)\frac{\partial \phi(x)}{\partial x^k}
F^2(x,\dx)
\right\} \,.
\label{eq:geod-coeff-conf}
\end{eqnarray}
If we write the conformal factor as
$\overline{F}^2(x,\dx)=\psi (x)F^2(x,\dx)$
we obtain
\begin{eqnarray}
	\label{eq:Gconfpsi}
\overline{G}^i(x,\dx)=
G^i(x,\dx)+
\frac{1}{4\psi(x)}
\left\{
2\frac{\partial \psi (x)}{\partial x^l}\dx^l\dx^i
-
g^{ik}(x,\dx)\frac{\partial \psi(x)}{\partial x^k}
F^2(x,\dx)
\right\}\,.	
\end{eqnarray}
Let $s \mapsto \tx(s)$ be a geodesic of the Finsler metric
$\overline{F}$ with
$\overline{F}^2(\tx,\tx')=1=2(E-U(x)) F^2(\tx,\tx'),$
i.e. 
\begin{equation*}
(\tx^i)''+2\overline{G}^i(\tx,\tx')=0\,.	
\end{equation*}
Choose 
\begin{equation*}
t(s)=t(s_0)+
\frac{1}{2}\int_{s_0}^s \frac{d\sigma}{E-U(\tx(\sigma))}\,,	
\end{equation*}
hence
$$
\frac{dt}{ds}=\frac{1}{2(E-U(\tx(s))}\,.
$$
For the inverse function $t \mapsto s=s(t)$ and the reparametrized
curve $x(t)=\tx(s)=\tx(s(t))$ we obtain
\begin{eqnarray*}
\dot{s}=\frac{ds}{dt}=2(E-U(x(t)))\,;\,
\ddot{s}=-2 \dot{U}(x)=-2 \frac{d}{dt} U(x(t))	\,,
\end{eqnarray*}
then $\dot{x}=\dot{x}(t)=
\tx'(s)\dot{s}=2(E-U(x)) \tx',$
i.e.
$F^2(x,\dx)=4(E(E-U(x))^2 F^2(\tx,\tx')=
2(E-U(x))\,.$
We compute the second derivative
using Equation~\eqref{eq:Gconfpsi}:
\begin{eqnarray*}
\ddx^i= (\tx^i)''\dot{s}^2+
(\tx^i)' \ddot{s}
=4(E-U(x))^2(\tx^i)''-2 \dot{U} (\tx^i)'\\
=-8 (E-U(x))^2 \overline{G}^i(\tx,\tx')
-\dot{U} \frac{\dx^i}{E-U(x)}	
=-2 \overline{G}^i(x,\dx)-\dot{U}\frac{\dx^i}{E-U(x)}\\
=
-2 G^i(x,\dx)+\frac{1}{E-U(x)}
\left\{
 \dot{U}\dx^i
 -  g^{ik}(x,\dx) \frac{\partial U}{\partial x^k} (E-U(x))
\right\}
-\dot{U}\frac{\dx^i}{E-U(x)}\\
= 
-2 G^i(x,\dx)-g^{ik}(x,\dx) 
\frac{\partial U}{\partial x^k}\,.
\end{eqnarray*}
We conclude from Equation~\eqref{eq:orbit-lagrange}
that
the reparametrized curve $t \mapsto x(t)$
is a solution of the Lagrangian system
defined
by $L(x,\dx)=F^2(x,\dx)/2-U(x).$

Finally we obtain the following statement, which is
well known in the case of a Riemannian metric,
cf. for example~\cite[Sec.4.1, Sec.4.2]{AKN}.
\begin{proposition}[Jacobi-Finsler metric]
	\label{pro:jacobi-finsler}
	Let $M$ be a differentiable manifold with a Finsler metric
	$F$ and a smooth function $U:M \longrightarrow \R.$
	Then a curve $t \mapsto x=x(t)$ resp. $
	t\mapsto (x(t),\dx(t))$ is a solution of the
	Lagrange equation for the Lagrangian
	$L(x,\dx)=F^2(x,\dx)/2-U(x)$
	with total energy $E$
	if and only if
	the reparametrized curve $s \mapsto \tx(s)=x(t(s))$
	with $(E-U(\tx(s))F^2(\tx,\tx')=1$ for all $s$
	is a geodesic of the
	\emph{Jacobi Finsler metric} 
	$\overline{F}^2(x,\dx)=2(E-U(x))F^2(x,\dx).$
\end{proposition}
\begin{remark}
	A similar statement also holds for the orbits of a
	Hamiltonian system of classical type. 
	The correspondence between the different descriptions is
	given by the following formulas for the 
	Lagrangian $L$, Hamiltonian $H$, the Finsler-Jacobi metric
	$\overline{F}$ and the Legendre transform
	$\mathcal{L}(x,\dx)=(x,y):$
	\begin{eqnarray}
	L(x,\dx)=\frac{1}{2}g_{ij}(x,\dx)\dx^i\dx^j-U(x)\,;\,
	H(x,y)=\frac{1}{2} g^{ij}(x,y)y_iy_j+U(x);\\
	y_i=g_{ij}(x,\dx)\dx^j\,;\,
	\dx^j=g^{jk}(x,y)y_k;\\
	\overline{F}^2(x,\dx)=2(E-U(x))F^2(x,\dx)
 =2(E-U(x))g_{ij}(x,\dx)\dx^i\dx^j		\,.
	\end{eqnarray}
	If the kinetic term is given by a Riemannian metric
	(i.e. $g_{ij}(x,\dx)=g_{ij}(x), g^{ij}(x,y)=g^{ij}(x)$) these formulas
	are standard, 
	see for example~\cite[Sec.3]{We} and
	\cite[Rem.2.7]{Co2022}.
	These formulae also hold in the Finsler case mainly due to
	the symmetries of the Cartan tensor,
	see Equation~\eqref{eq:cartan-symmetry},
	this was already shown by Ruiz~\cite[Sec.6]{Ru}.
\end{remark}
\begin{example}
	\label{example:harmonic-oscillator}
	The \emph{harmonic oscillator} on $\R^{2n}
	=\{(x^1,\ldots,x^n,\dx^1,\ldots,\dx^n);
	x^i, \dx^j\in \R\}$ with frequencies
	$\alpha_i>0,i=1,\ldots,n$  is defined by the
	Hamiltonian
	$$
	H(x,y)= \frac{1}{2}\sum_{i=1}^n (y_i^2 + \alpha_i^2(x^i)^2)
	$$
	resp. the Lagrangian
	$$
	L(x,\dx)=\frac{1}{2}\sum_{i=1}^n \left\{(\dx^i)^2
	-\alpha_i^2 (x^i)^2\right\}
	$$ with Euclidean norm $F^2(x,\dx)=\sum_i (\dx^i)^2$
	and 
	potential $U(x)=(\sum_i \alpha_i^2(\dx^ i)^2)/2.$
	Then the Lagrange equations~\eqref{eq:lagrange} 
		imply
	$\ddx^i+\alpha^2_i x^i=0$ for all $i=1,\ldots,n.$
	In particular, if the quotients $\alpha_i/\alpha_j, i\not=j$ are all
	irrational
	(called 
	\emph{weakly non-resonant harmonic oscillator}) there are exactly $n$ 
	geometrically distinct and prime periodic orbits for fixed energy
	$E>0$
	with minimal periods $T_i=2\pi/\alpha_i.$ 
	These periodic orbits $x_j=x_j(t), j=1,\ldots,n$ 
	with $x_j(0)=0$ are given by
	\begin{equation}
		\label{eq:xj}
		x_j^j(t)=\frac{\sqrt{2E}}{\alpha_j}\sin(\alpha_j t), x_j^i(t)=0, i\not=j.
	\end{equation}
	These periodic orbits are brake orbits with
	brake points resp. rest points $x^j_j(\pm\pi/(2\alpha_j)) =\pm \sqrt{2E}/\alpha_j,
	x^i_j(\pm \pi/(2\alpha_j))=0, i\not=j.$
	A Jacobi field $v=v(t)=v^i \partial_i$ along $x_j(t)$ orthogonal to $\dx_j$
	with respect to the Euclidean metric determined by $F$
	solves the Jacobi Equation
	\begin{equation*}
		\ddot{v}^i(t) + \alpha_i^2 v^i(t)=0	
	\end{equation*}
	(no summation convention here),
	i.e.
	$v^i(t)=v^i(0) \cos(\alpha_i t)+\dot{v}^i(0)\sin(\alpha_i t).$
	If $\alpha_i/\alpha_j$ is irrational for all $i\not=j$ then
	there is no periodic Jacobi field $v=v(t)$ along $x_j$ orthogonal
	to $x_j,$ in particular the brake orbits and
	all their iterates are non-degenerate, the linearized Poincar\'e
	mapping of the orbit $x_j$ is a linear symplectic mapping which is the
	sum of $(n-1)$ two-dimensional rotations with
	angles $\alpha_i, i\not=j.$
	The brake orbits $x_j,  j=1,\ldots,n$ (cf. Equation~\eqref{eq:xj}) 
	are simple, i.e. the restriction
	$x_j:[-\pi/(2\alpha_j),\pi/(2\alpha_j)]\longrightarrow 
	\{x \in \R^n, \sum_i \alpha_i^2 (x^i)^2\le 2 E\}$
	is injective and all $n$ brake orbits have the origin $0 \in \R^n$
	as common intersection point.	
	For a \emph{resonant harmonic oscillator} there are two
	indices $i,j$ (we set without loss of generality
	$i=1,j=2$) and $\alpha>0$ with
	$\alpha_1= m_1 \alpha, 
	\alpha_2=m_2 \alpha$ for relatively prime
	integers $m_1,m_2.$ 
	
	Then we get as a particular brake orbit $x_*$
	with $a^1,a^2 \in \R:$
	$$
	x_*^1(t)= a^1 \cos(m_1\alpha t);
	x_*^2(t)=a^2 \cos(m_2 \alpha t), x_*^k(t)=0, k\ge 3.
	$$
	The total energy of $x_*$ 
	is given by
	$E=\alpha_1 (a^1)^2+\alpha_2(a^2)^2.$
	Hence the minimal period is $2\pi/\alpha$
	and this brake orbit represents a 
	\emph{Lissajous figure,}
	in particular this prime brake orbit is not simple.
	In addition one can check that the orbit is 
	degenerate:
	There is a one-parameter family of brake orbits
	$x_{*,s}$
	of energy $E=\alpha_1 (a^1)^2+\alpha_2(a^2)^2$ 
	and minimal period $2\pi/\alpha$
	with
	\begin{eqnarray*}
		x_{*,s}^1(t)=  a^1\{\cos (s) 
		\cos(m_1\alpha t)+ \sin (s) \sin(m_1\alpha t)\};\\
		x_{*,s}^2(t)=a^2 \cos(m_2 \alpha t), x_{*,s}^k(t)=0, k\ge 3,
	\end{eqnarray*}
	which shows that the brake orbit
	$x_{*,0}=x_*$ is degenerate.
\end{example}
\section{Perturbing a single segment of an orbit}
Let $v,w\in T_pM$ be a unit vectors 
for which $g(v)(v,w)=\mathcal{L}(v)(w)=0.$
Here $\mathcal{L}$ is the Legendre transformation introduced
in Equation~\eqref{eq:Legendre}.
Let $(x^1,\ldots,x^n)$ be coordinates defined for
$x^1 \in [-2\eta,2\eta], (x^2)^2+\ldots+(x^n)^2<\epsilon$
and 
$g_{ij}(x^1,0,\ldots,0,1,0,\ldots,0)=\delta_{ij}$
such that the curves $t \mapsto (t,x^2,\ldots,x^n)$ are geodesics,
and $s \mapsto s (0,x^2,\ldots,x^n)$ are geodesics, too.
The geodesics $t \mapsto (t,x^2,\ldots,x^n)$ can be chosen to
be orthogonal to the hypersurface $t=0$ with respect
to the 	Riemannian metric $g_{ij}((0,x^2,\ldots,x^n),
(1,0,\ldots,0)).$ Note that in contrast to the Riemannian case in the Finsler case 
orthogonality depends also on the direction.
And the coordinate line $\gamma(t)=\gamma_0(t)=(t,0,\ldots,0)$
is the geodesic with $\gamma(0)=p=(0,\ldots,0),
\gamma'(0)=(1,0,\ldots,0)=v,$ and 
$s \mapsto \delta(s)=(0,s,0,\ldots,0)$ is the geodesic 
with $\delta(0)=p$ and $ \delta'(0),w$ are linearly dependent.
In the Riemannian case they can be chosen to be equal, in the
Finsler case $g^v(w,w)=1$ which does not imply
$F^2(w)=g^w(w,w)/2=1.$
We denote by $\xi=\xi_{v,w}$ the diffeomorphism
defined by these coordinates,
i.e. $\xi^{-1}(q)=(x^1,\ldots,x^n)$ as a chart
parametrizing the tubular neighborhood
$Tb_v(\eta.\epsilon)$ of the
geodesic $\gamma=\gamma(t):$
\begin{equation}
	\label{eq:xi}
\xi=\xi_{v,w}: 
[-2\eta,2\eta]\times D^{n-1}(\epsilon)
\longrightarrow Tb_v(\eta,\epsilon)\,;\,
\xi^{-1}(q)=(x^1(q),\ldots,x^n(q))\,.
\end{equation}
Let $\alpha:[2\eta,2\eta]\longrightarrow [0,1]$
be a smooth cut-off function with $\alpha(t)=0$ for
$|t|\ge \eta+4\epsilon$ and $\alpha(t)=1$ for
$|t|\le \eta +2\epsilon.$
Assume that the coordinates $x=(x^1,\ldots,x^n)$
are defined on $x^1\in [-2\eta,2\eta], 
\sum_{i=2}^{n}
(x^i)^2< \epsilon^2.$
Then for sufficiently small $\delta>0$ there 
is a one parameter group of diffeomorphisms
\begin{equation}
\tilde{\Psi^s}: (x^1,\dots,x^n)
\in [-2\eta,2\eta]\times D^{n-1}(\epsilon)
\mapsto
(\tilde{x}^1,\ldots,\tilde{x}^n)
\in
[-2\eta,2\eta]\times D^{n-1}(\epsilon)
\end{equation}
such that
\begin{equation}
\tilde{\Psi}^s(t,0,\ldots,0)=
(t,u_s(t),0,\ldots,0)
\end{equation}
for an 
even function, i.e. $u_s(-t)=u_s(t),$ 
for which the restriction 
$u_s:[0,2\eta]\longrightarrow [0,s]$
is monotone decreasing and satisfies
$u_s(t)=s$ for all $t\in [0,2\eta+2\eta]$
and $u_s(t)=0$ for all $t \ge \eta+4\epsilon.$
In addition for sufficiently small
$s \in (0,\delta]$ the 
diffeomorphism lies in an arbitrary small
neighborhood of the identity 
with respect to the $C^{\infty}$-topology.
For the details of the construction see
the Proof~\cite[Lemma 1]{R24}.
For the coordinates $(\tilde{x}^1,\ldots, \tilde{x}^n)$
(depending on $s$)
we define a new coordinate system
$z^i=z^i(\tilde{x}^1,\ldots,\tilde{x}^n)$
such that $z^i(t,0,\ldots,0)=\tilde{x}^i(t,0,\ldots,0)$
and such that for the corresponding 
metric coefficients
$$
g_{ij}(z,\dot{z})=
\frac{1}{2}\frac{\partial F^2(z,\dot{z})}{\partial \dot{z}^i
	\partial \dot{z}^j}
$$
we have
\begin{equation}
	\label{eq:gij-z}
g_{ij}((t,0,\ldots,0),(1,0,\ldots,0))=
\delta_{ij}.
\end{equation}
 It also follows that for sufficiently small
$\delta$ the diffeomorphism $z=z(\tilde{x})$ 
resp. $z^i=z^i(\tilde{x}^j)$ lies in an
arbitrary small neighborhood of the identity with
respect to the $C^{\infty}$-topology.

For our computations we now use the coordinates
$z=(z^1,\ldots,z^n),$ along the coordinate line
$\tgamma(t)=\gamma_s(t)=
(t,0,\ldots,0), -2\eta \le t \le 2\eta,$ 
which also depend on $s\in [0,\delta],$
and the metric along $\gamma$ satisfies
Equation~\eqref{eq:gij-z}.

In particular the coordinate line 
$\tgamma(t)$ is parametrized by arc length.
Then the geodesic coefficients of 
$\tgamma'(t)=((t,0,\ldots,0),(1,0,\ldots,0))$
are given by
\begin{eqnarray}	
G^i(\tgamma'(t))=\frac{1}{4}\delta^{il}\left\{
2\frac{\partial g_{jk}(\tgamma'(t))}{\partial z^l}-
\frac{\partial g_{jl}(\tgamma'(t))}{\partial z^k}
\right\}\delta^{j1}\delta^{l1}
\nonumber
\\
\label{eq:Gi}
=\frac{1}{4}\left\{
2\frac{\partial g_{1i}(\tgamma'(t))}{\partial z^1}-
\frac{\partial g_{11}(\tgamma'(t))}{\partial z^i}
\right\}
=- \frac{1}{2} 
\frac{\partial g_{11}(\tgamma'(t))}{\partial z^i}\,.
\end{eqnarray}
Hence 
\begin{equation}
	\label{eq:G1}
G^1(\tgamma'(t))=0.
\end{equation}
For the geodesic coefficents $\overline{G}^i(\gamma'(t))$
of the conformally equivalent Finsler metric
$\overline{F}(z,\dot{z})=\exp(\phi(z))F(z,\dot{z})$ we obtain
from Equation~\eqref{eq:geod-coeff-conf}:
\begin{eqnarray}
\overline{G}^i(\tgamma'(t))=
G^i(\tgamma'(t))+\frac{1}{4}\delta^{ik}\left\{
2 \frac{\partial \phi}{\partial z^l}\delta_{jk}
-\frac{\partial \phi}{\partial z^k}\delta_{jl}
\right\}\delta_{j1}\delta_{l1}
\nonumber
\\
\label{eq:Gi4}
=	G^i(\tgamma'(t))+\frac{1}{4}\exp(-\phi(\tgamma(t)))
\left\{2 \frac{\partial \phi(\tgamma(t))}{\partial z^1}
-\frac{\partial \phi(\tgamma(t))}{\partial z^i}
\right\}\,.
\end{eqnarray}
We are looking for a function $\phi,$ such that
$\tgamma$ is a geodesic for the metric $\overline{F},$
hence $\overline{G}^i(\tgamma'(t))=0$ for all $1\le i\le n.$

For $i=1$ we obtain from Equation~\eqref{eq:G1} 
- that $\phi$ is constant
along $\tgamma.$ Let $\phi(\tgamma(t))=0$
for all $t,$ then we obtain
from Equation~\eqref{eq:Gi} and Equation~\eqref{eq:Gi4}:
\begin{equation}
\frac{\partial \phi (\tgamma(t))}{\partial z^i}=
2 \frac{\partial g_{11}(\tgamma'(t))}{\partial z^i}\,.    
\end{equation}
If $\alpha:\R\longrightarrow [0,1]$
is a smooth function with $\alpha(t)=0$
for $|t|\ge \epsilon$ and $\alpha(t)=1$
for $|t|\le \epsilon/2$ choose
for $r^2=(z^1)^2+\ldots+(z^n)^2:$
\begin{equation}
\phi(z^1,\ldots,z^n))=
\alpha(r^2)
\frac{\partial g_{11}(\tgamma'(t))}{\partial z^i}z^i
\end{equation}
Hence for the Finsler metric $\overline{F}$
the curve $\gamma(t)$ is a geodesic and $\overline{F}$
coincides with the original metric $F$ outside a small 
tubular neighborhood of the segments $\gamma|[\eta,\eta+4\epsilon]$
and $\gamma|[-\eta-4\epsilon,-\eta].$
Hence we obtain the following
\begin{lemma}[Perturbation Lemma]
		\label{lem:perturbation-conformal}
	Let $(M,F)$ be a Finsler metric, $p\in M,v,w \in T_pM$
	be two unit vectors orthogonal to each other with respect
	to $g^v.$ Then for sufficently small $\eta, \epsilon>0, 
	7\epsilon <\eta$  such that the
	map $\xi=\xi_{v,w}$  defined in Equation~\eqref{eq:xi}
	is a diffeomorphism and the coordinate line
	$\gamma(t)=(t,0,\ldots,0), -2\eta\le t\le 2\eta$ 
 is a geodesic parametrized 	by arc length, we have the following: 
	
	There is a smooth one-parameter family of 
	smooth functions $\phi=\phi_s$ with $\phi_0=0$
	such that for the 
	conformally equivalent metric
	$\overline{F}=\exp(\phi) F$ 
	the curve
	$\gamma_s:[-2\eta,2 \eta]\longrightarrow Tb_v(\eta,\epsilon)$
	with $\gamma_s(t)=\xi_{v,w}(t,u_s(t),0,\ldots,0)$ and 
	$u_s(t)=s$ for $|t|\le \eta+2\epsilon$ and
	$u_s(t)=0$ for $|t|>\eta+4\epsilon$ 
	is a geodesic parametrized by arc length.
	And the support of $\phi$ lies in 
	an arbitrary small neighborhood $U$ of 
	$\gamma([-\eta-4\epsilon,-\eta]) \cup
	\gamma([\eta,
	\eta+4\epsilon]).$
	Given a neighborhood of the $0$-function in the space
	of smooth functions with the 
	$C^{\infty}$-topology one can choose
	$s$ sufficiently small such that one can find
	$\phi$ lying in this neighborhood with support in $U.$	
\end{lemma}
\begin{remark}
	The proof also shows that by a conformal change of the
	Finsler metric an arbitrary smooth curve 
 parametrized by arc length and
 	without self-intersections
	can be made a geodesic curve.	
\end{remark}
\section{Intersection of orbits of Lagrangian systems}
\label{sec:intersection}
We investigate self-intersection points of orbits of
Lagrangian systems and the intersection of distinct
orbits and obtain results similar to the corresponding
results for closed geodesics of Riemannian metrics and
reversible Finsler metrics, cf.~\cite[Sect.3]{R24}
resp. for geodesic segments,
cf.~\cite[Sec.2]{R24a}.
\begin{lemma}
\label{lem:intersection}
Assume that a Lagrangian system of classical type
with reversible Finsler metric $F$ and potential $U$ is
given on a manifold $M$. 
Let $E$ be regular value of the potential $U$ with 
compact potential well $M^ {\le E}=\{x\in M; U(x)\le E\}.$
Let $\gamma: \Stau\longrightarrow M^{\le E},
\gamma_*: S^1_{\tau_*}\longrightarrow M^{\le E}$ 
be geometrically distinct
and prime periodic orbits of
energy $E$ of period $\tau$ resp. $\tau_*.$
\begin{itemize}
\item[(a)] The set of \emph{self-intersection points}
$$
DP(\gamma)=\#\{\gamma(t), \exists s\in \Stau, s\not= t,
\gamma(t)=\gamma(s), \gamma'(t)\not=-\gamma'(s)\}\,.
$$
is finite.
\item[(b)] The intersection $I(\gamma,\gamma_*)=
\gamma(\Stau)\cap \gamma_*(S^1_{\tau_*})$ is finite.	
\end{itemize}
\end{lemma}
The Lemma follows from 
\cite[Lemma 2]{R24} resp.
\cite[Lemma 1]{R24a} 
by the following Remark showing the relation
of brake orbits with orthogonal geodesic chords:
\begin{remark}
\label{rem:ogc}
We obtain from~\cite[Theorem 1.9]{Co2022}
the following (the Riemannian case was shown
in~\cite{GGP2004}):
Using the Jacobi-Finsler metric 
$\overline{F}^2=2(E-U)F^2$
one can define a function
$\psi: M^{\le E}\longrightarrow \R$
as follows: For $y \in M^{\le E}$ let $X_y$ be the set
of curves $\sigma\in H^1([0,1],M^{\le E})$ with
$y=\sigma(0), \sigma(1)\in M^E.$
Then we define the function
$$
\mathcal{T}(\sigma)=\frac{1}{2}\int_0^1
\overline{F}^2(\sigma(s),\sigma'(s))\,ds.
$$
Then the function $\psi: M^{\le E}\longrightarrow \R,
\psi(y)=\inf\{\mathcal{T}(\sigma); \sigma \in X_y\}$
is smooth near the boundary $M^{E}=\psi^{-1}(0)$ 
and $0$ is a regular value.
Given a brake orbit $\gamma:\Stau \longrightarrow M^{\le E}$ of
energy $E$ with $U(\gamma(0))=U(\gamma(\tau/2))=E$
there is a sufficiently small $\delta>0$ and a subinterval
$[\alpha,\beta]\subset [0,T/2]$ 
such that $\psi(\gamma(\alpha))=\psi(\gamma(\beta))=\delta$
and $\psi(\gamma(t))> \delta$ for all $t\in (\alpha,\beta).$
And $\gamma|[\alpha,\beta]$ is (up to parametrization) 
an \emph{orthogonal geodesic chord} of 
$M':=\{x\in M\,;\,\psi(x)\ge \delta\}.$
Hence for $t=0,T/2$ the tangent space
$T_{\gamma(t)}\psi^{-1}(\delta)$ equals the kernel
of $\mathcal{L}(\gamma'(t)).$ Here $\mathcal{L}$ is the
Legendre transformation introduced in Equation~\eqref{eq:Legendre}.
The Jacobi-Finsler metric $\overline{F}^2=2(E-U(x))F^2$ is 
well defined on $M',$ the injectivity radius is positive
and hence the arguments of 
\cite[Lemma 2]{R24} resp. \cite[Lemma 1]{R24a} give the proof
of Lemma~\ref{lem:intersection}.
\end{remark}
\section{Proof of Theorem~\ref{thm:main}}
	\label{sec:proof}
It follows from~\cite[Theorem 1]{Be} that for a generic
function $U \in C^{\infty}(M)$ the energy level
$H^{-1}(E)$
of $L=F^2-U$ with $H=F^2/2+U$ is regular and all periodic orbits 
of energy $E$ are non-degenerate.
It also follows that for some $a>0$ the set $\mathcal{U}(a)$
of potentials
$U \in C^{\infty}(M)$ for which all periodic orbits of
period $\le a$ of energy $E$ are non-degenerate,
is an open and dense subset with respect to the
$C^{\infty}$-topology. Since $M$ (and hence $M^{\le E}$)
is compact, and all periodic orbits of energy $E$ and 
period $\le a$ are non-degenerate we conclude that there are 
only finitely many periodic orbits
of period $\le a$ and energy $E.$ 
Let $\mathcal{U}^*(a)$ be the set of smooth functions
$u \in C^{\infty}(M)$ such that all prime periodic orbits
in $H^{-1}(E)$
do not have self-intersections and do not intersect each other.
Since there are only finitely many geometrically distinct
periodic orbits of period $\le a$ 
in $H^{-1}(E)$ the set
$\mathcal{U}^*(a)$ is an open subset of $C^{\infty}(M).$
Let $\gamma_1,\ldots,\gamma_r$ be the periodic orbits of
energy $E$ and period $\le a$
up to geometric equivalence. Hence these are 
up to parametrization either closed
geodesics of the Jacobi-Finsler metric $\overline{F}$ or
brake orbits, cf. Remark~\ref{rem:ogc}.
We still have to show that $\mathcal{U}^*(a)$ is a dense
subset of $C^{\infty}(M):$
We conclude from Lemma~\ref{lem:intersection} that the union
$$DP(a)=DP(\gamma_1)\cup \ldots \cup DP(\gamma_r)
\cup
\bigcup_{j\not=k} I(\gamma_j,\gamma_k)
$$
of double points $DP(\gamma_j),j=1,\ldots,r$
of the prime periodic orbits of energy $E$ and period $\le a$
and the intersection points $I(\gamma_j,\gamma_k)$ for
distinct $j,k$ is a finite set. Now we consider the orbits
$\gamma_j, j=1,\ldots,r$ after a reparametrization as 
geodesics of the Jacobi-Finsler metric as 
described in Remark~\ref{rem:ogc}.
Using the function $\psi:M^{\le E}\longrightarrow \R$
for sufficiently small $\delta$ the periodic orbits of
energy $E$ and period $\le a$ correspond to either
closed geodesics 
or orthogonal geodesic chords
in $M'=
\{x\in M^ {\le 0}; \psi(x) \ge \delta\}.$

For sufficiently small
$\eta>0$ and for a point $p \in DP(T)$ the intersection
$B_p(2\eta)\cap\bigcup_{j=1}^r \gamma_j(S^1_{T_j})$
consists with respect to the Jacobi-Finsler metric of
$N$ geodesic segments $\gamma_j, j=1,\ldots,N$ with common midpoint 
$p=\gamma_j(0)$ for 
some $N \in \mathbb{N}.$
For $0 <7\epsilon<\eta$ we find open neighborhoods
$U^{-}_j$ of $\gamma_j([-\eta-4\epsilon,-\eta])$
resp. $U^{+}_j$ of $\gamma_j([\eta+2\epsilon,\eta+4\epsilon])$
which are all pairwise distinct
such that the following holds by Lemma~\ref{lem:perturbation-conformal}:
In any neigbhorhood of the zero function in 
$C^{\infty}(M)$ we find a function $\phi$ with
support in the union 
$U^*=\bigcup_{j=1}^N U_j^{\pm}$ of the sets $U_j^{\pm},j=1,\ldots,N$ 
such that for the conformally equivalent metric 
$ \exp(\phi)\overline{F}^2$ the geodesics $\gamma_j,j=1,\ldots,r$
are perturbed to geodesics $\tilde{\gamma}_j$
with the same length and without
self-intersections and without intersections between each other.
Then these geodesics are (up to parametrization) solutions
of the reversible Lagrangian system defined by
the Finsler metric $F$ and the potential
$(1-\exp(\phi))E +\exp(\phi)U.$ 
This follows from:
$\exp(\phi)\overline{F}^2=2\exp(\phi))(E-U)F^2
=2(E-\{(1-\exp(\phi))E+\exp(\phi)U\}F^2.$

Hence we have shown that the set $\mathcal{U}^*(a)$ is dense in $C^{\infty}(M).$
\section{Proof of Corollary~\ref{cor:periodic}}
\label{sec:periodic}
Let $\mathcal{U}_1$ be an open and dense subset of potentials 
in $C^{\infty}(M)$ 
for which
$E$ is a regular value of $U,$ 
and $E$ is a regular value of $H=\kin+U.$
For $U \in \mathcal{U}_1$ it follows from~\cite[Thm.A]{GZ}
that there is a periodic orbit $\gamma$
in $H^{-1}(E),$ which is a brake orbit if
$M^E$ is non-empty.
If $\tau_1$ is the minimal period of this periodic orbit
$\gamma$
we define the following subset $\mathcal{U}\subset 
\mathcal{U}_1:$ For $U \in \mathcal{U}$ all
periodic orbits in $(\kin+U)^{-1}(E)$
of minimal period $\le 2 \tau_1$ are non-degenerate
and simple.
It follows from Theorem~\ref{thm:main} that the set
$\mathcal{U}$ is an open and dense subset
of $\mathcal{U}_1.$


\begin{thebibliography}{999999}
\bibitem{AS}
Albu-Sch\"affer, Alin, and Sachtler, Arne,
What can algebraic topology and differential geometry teach us
about intrinsic dynamics and global behaviour of robots?
In:\emph{Robotics research,} edited by
Billard, Aude, Asfour, Tamim, and Khatib, Oussama,
Springer Proc.~Adv.~Robotics (SPAR. volume 27), 468--484
Springer Nature, 2023
\bibitem{AKN}
Arnold, Vladimir I., Kozlov, Valery V., and 
Neishtadt, Anatoly I.,
\emph{Mathematical aspects of classical and
	celestial mechanics,}
transl. from the Russian by E.Khukhro,
3rd revised ed., Encycl.~Math.~Sciences, Vol. 3, Dyn.~Systems III,
Springer 2006
\bibitem{Be} Bernard, Patrick,
Non-degeneracy of closed orbits for generic
potentials, \emph{J.~\'Ecole~Polytechn.~Math.} 11 (2024) 363-393
\bibitem{BK} Bolotin, S.V., and Kozlov, V.V.,
Libration in systems with many degrees of freedom,
\emph{J. Appl. Math. Mech.} 42, 256--261 (1979); 
translation from \emph{Prikl. Mat. Mekh.} 42, 245--250 (1978) (Russian)
\bibitem{Co2022}
Corona, Dario, and Giannoni, Fabio,
Brake orbits for Hamiltonian systems of the classical type
via geodesics in singulare Finsler metrics,
\emph{Adv.~Nonlin.~Anal.} 11 (2022) 1223--1248;
\bibitem{DdP}
Dahinden, Lucas, and de Pooter, Jacobus,
On the intersections of projected Hamiltonian orbits
in cotangent bundels,
{\tt arxiv:2602.15693v1} (2026)
\bibitem{GGP2004}
Giamb\`o, Roberto, Giannoni, Fabio, and Piccione, Paolo,
Orthogonal geodesic chords, brake orbits and homoclinic orbits
in Riemannian manifolds,
\emph{Adv.~Differ.~Equat.} 10 (2004) 1--24
\bibitem{GGP2022}
Giamb\`o, Roberto, Giannoni, Fabio, and Piccione, Paolo,
Multiple orthogonal geodesic chords and a proof of
Seifert's conjecture on brake orbits, 
{\tt arXiv:2002.09687v3}(2022)
\bibitem{GZ}
Gluck, Herman and Ziller, Wolfgang,
Existence of periodic motions of conservative systems,
In: \emph{Semin. on minimal submanifolds,} 
{Ann}. {Math}. {Stud}. 103 (1983) 65-98
\bibitem{Kn} 
Knauf, Andreas,
Closed orbits and converse KAM theory,
\emph{Nonlinearity} 3 (1990) 961--973
\bibitem{Ko}
Kozlov, V. V.,
Calculus of variations in the large and classical mechanics,
\emph{Russ.~Math.~Surv.} 40 (1985) 37--71;
translation from
\emph{Usp.~Mat.~Nauk} 40
(1985) 33--60 (Russian)
\bibitem{Ma}
Ma\~n\'e, R.: Generic properties and problems of 
minimizing measures of Lagrangian systems.
\emph{Nonlinearity} 9 (1996) 273-310
\bibitem{Ol} Oliveria, F.R.,
Generic properties of Lagrangians of surfaces: the 
Kupka-Smale theorem,
\emph{Discrete Contin. Dynam. Systems} 21 (2008) 551--569
\bibitem{R24} Rademacher, Hans-Bert,
Simple closed geodesics in dimensions $\ge 3,$
\emph{J.Fixed Point Theory Appl.},\textbf{26,} 5 (2024)
\bibitem{R24a}
Rademacher, Hans-Bert, Geodesic loops and orthogonal geodesic chords without self-intersections, 
\emph{Nonlinear~Anal.} TMA 263 (2026) 113952
\bibitem{RR} Rifford, L., and Ruggiero, R.O.,
Generic properties of closed orbits of Hamiltonian
flows from Man\*e's viewpoint, \emph{Intern.~Math.~Res. Notices}
(2012) 5246--5265
\bibitem{Ru} Ruiz M., Otto Raul,
Existence of brake orbits in Finsler mechanical systems.
\emph{Geom.Topol.,III. Lat.~Am.~Sch.~Math.~Proc.,
	Rio de Janeiro 1976, Lect.~Notes~Math.} 597 (1977) 542--567
\bibitem{Se} Seifert, Herbert, 
Periodische Bewegungen mechanischer Systeme,
\emph{
Math.~Zeitschr.} 51 (1948) 197--216 
\bibitem{Sh} Shen, Zhongmin,
\emph{Lectures on Finsler geometry,}
World Scientific, Singapore, 2001
\bibitem{We} Weinstein, A., Periodic orbits for convex hamiltonian systems,
\emph{Ann.~Math.} 108 (1978) 507--518
\end{thebibliography}
\end{document}